\documentclass[12pt]{article}
\usepackage{amssymb}
\def \Z{{\mathbb Z}}
\def \p{{\mathbb P}}

\begin{document}
\title{On morphisms onto quadrics}
\author{Ekaterina Amerik}
\date{}
\maketitle

The purpose of this note is to prove the following result:

\medskip

{\bf Theorem} {\it Let $X$ be a smooth projective variety of dimension
$n$ with cyclic N\'eron-Severi group, and let $Q$ denote a smooth $n$-dimensional 
quadric ($n\geq 3$).
The degree of a morphism $f:X\rightarrow Q$ is bounded in terms of the
discrete invariants of $X$. 

More precisely, let $L$ denote an ample generator of $NS(X)$, and let 
$k, r$ be such that $K_X\equiv kL$ and $rL$ is very ample. If 
$f^*{\cal O}_Q(1)\equiv mL$, then $m<2k+(2n+5)r$ when $n\geq 4$, and 
$m<3k+16r$ for $n=3$.}

\medskip

Obviously, the degree of $f$ is equal to $m^nL^n/2$.

\medskip

The existence of such a result was briefly mentioned in \cite{A}, Remark 4.1.2;
however, the proof (very much analogous to that of \cite{A}, Theorem
4.1.1) was never 
made explicit. Recently, J.-M. Hwang, to whom I am grateful
for his interest in that work,
suggested that I should write down this proof, which I am glad to do
here.

\medskip

Note that Hwang and Mok \cite{HM} proved the boundedness of the
degree
of morphisms between Fano manifolds with cyclic Picard group
admitting rational curves with trivial normal bundle. The quadric, of
course,
does not have such curves: its anticanonical class is ``too
positive''.
More generally, the method of Hwang and Mok does not seem to be
applicable when the target is a homogeneous variety. 

The theorem above can be understood as an indication that the degree
of a morphism $f:X\rightarrow Y$, where $X,Y$ are $n$-dimensional
varieties with Picard number one and $Y$ is Fano, is bounded
unless
$Y=\p^n$ (compare with \cite{ARV}, Proposition 2.1, treating the
case 
when $Y$ has trivial
canonical class). However, I do not know whether it is possible to
generalize the argument below, say, to some other homogeneous varieties. 

\

Let $f:X\rightarrow Q$, $NS(X)\cong \Z$, be a morphism; $f$ is 
finite if non-constant. Take $l$ a sufficiently general line on $Q$, 
and let $S\subset Q$ be a
general 2-dimensional linear section of $Q\subset \p^{n+1}$ containing $l$ (so $S$ is a 
two-dimensional quadric). Write $S=Q\cap H_1\cap H_2 \cap \dots \cap H_{n-2}$,
with $H_i$ sufficiently general hyperplanes.

One has: $N_{l,S}={\cal O}_{\p^1}$ and 
$N_{S,Q}|_l={\cal O}_{\p^1}(1)^{\oplus n-2}$
That is, the exact sequence of normal bundles 
$$0\rightarrow N_{l,S}\rightarrow N_{l,Q}\rightarrow N_{S,Q}|_l\rightarrow 0$$
splits. Let $C=f^{-1}l$ and $M=f^{-1}S$. As $l$ and $S$ are chosen to be 
sufficiently general, $C$ and $M$ are smooth by Kleiman's version
of Bertini theorem (\cite{H}, chapter III, section 10). The exact sequence of normal bundles
$$0\rightarrow N_{C,M}\rightarrow N_{C,X}\rightarrow N_{M,X}|_C\rightarrow 0$$
is just the pull-back of the corresponding sequence on $Q$, that is, it
splits. So the induced map $\alpha:H^0(C, N_{C,X})\to H^0(C, N_{M,X}|_C)$
is surjective. 

In particular, $Im(\alpha)$ contains the image of
the restriction map $\beta: H^0(M, N_{M,X})\to H^0(C, N_{M,X}|_C)$. This is
interpreted as follows: any infinitesimal deformation of $M$ in $X$ contains
an infinitesimal deformation of $C$. That is, the class of $C$ in 
$H^{1,1}(M)$ is ``infinitesimally fixed'' according to the terminology of
\cite{CGGH}. 

Now there are various results of ``infinitesimal Noether-Lefschetz'' type,
which state that if $M$ is sufficiently ample on $X$, then all 
infinitesimally fixed classes in $H^{1,1}(M)$ are restrictions of (1,1)-classes
on $X$. Our class $[C]$ obviously cannot be such a restriction: indeed,
$N_{C,M}={\cal O}_C$, that is, $[C]^2=0$; but, thanks to the hypothesis
that $NS(X)=\Z$, any (1,1)-class which has zero square on $M$ and is a 
restriction
of a (1,1)-class on $X$, must be numerically trivial, which is obviously
not the case for $[C]$. So such results give the non-existence of $f$ 
as soon as $M$ is sufficiently ample, that is, $deg(f)$ is sufficiently
high.

In order to make this more precise, we shall prove the following proposition,
which is a very direct generalization of Proposition 3.4 from
\cite{EL}:

\medskip

{\bf Proposition} {\it Let $X$ be a smooth projective $n$-fold and
$Y$ a smooth complete intersection of ample divisors $D_1, D_2, \dots, D_{n-2}$
on $X$ (so that dim(Y)=2). Suppose that for some very ample $A$, one has
$D_i\in |K_X+(n+1)A+B_i|$
and $D=D_1+D_2+\dots + D_{n-2}\in |nK_X+(n+1)^2A+E|$ with $B_i$ and $E$ nef. 
Then the infinitesimal Noether-Lefschetz theorem holds for $Y\subset X$, that 
is, any (1,1)-class on $Y$ which
remains of type (1,1) under all infinitesimal deformations of $Y$ in $X$, 
is a restriction
of a class on $X$.}

\medskip

{\it Proof:} We follow the exposition of \cite{EL}. Let 
$N=\oplus_i{\cal O}_Y(D_i)$ be the normal bundle of $Y$ in $X$, so
that
$H^0(Y,N)$ is the space of infinitesimal deformations of $Y$ in $X$.
According to the theory of
\cite{CGGH}, the space of infinitesimally fixed (1,1)-classes is the
kernel of the
map

$$\alpha: H^1(Y, \Omega^1_Y)\stackrel{\gamma}{\longrightarrow}
H^2(Y, N^*)
\stackrel{\beta}{\longrightarrow} Hom(H^0(Y,N), 
H^2(Y, {\cal O}_Y)),$$
$\gamma$ being the usual coboundary map, and
 $\beta$ coming from the multiplication map 
$\mu:H^0(Y, N)\otimes H^0(Y, K_Y)\rightarrow H^0(Y, N\otimes K_Y)$ by Serre
duality. We shall be done if we check that

1) $\beta$ is injective, that is, $\mu$ is surjective;

2) the restriction map $H^1(X, \Omega^1_X)\rightarrow H^1(Y, \Omega^1_X|_Y)$
is surjective.

To check the surjectivity of $\mu$, it is enough to check the surjectivity
of each ``component'' 
$$\mu_i: H^0(Y, {\cal O}_Y(D_i))\otimes H^0(Y, K_Y) \rightarrow 
H^0(Y,{\cal O}_Y(D_i)\otimes K_Y).$$

From the resolution
$$(\ast)\  0\rightarrow {\cal O}_X(-D)\rightarrow \dots \rightarrow 
\oplus_{i<j} {\cal O}_X(-D_i-D_j)\rightarrow \oplus_i{\cal O}_X(-D_i)
\rightarrow {\cal O}_X \rightarrow {\cal O}_Y\rightarrow 0$$
one sees that
$H^0(X, {\cal O}_X(D_i)\otimes K_X\otimes {\cal O}_X(D))$ surjects
onto $H^0(Y,{\cal O}(D_i)\otimes K_Y)$ (indeed, the $D_i$ are ample,
so this is implied by the Kodaira vanishing).
That is, the surjectivity of $\mu_i$ would, by restriction, 
follow from the surjectivity
of the multiplication map on $X$:
$$H^0(X, {\cal O}_X(D_i))\otimes H^0(X, {\cal O}_X(D)\otimes K_X)
\rightarrow H^0(X, {\cal O}_X(D)\otimes K_X\otimes{\cal O}_X(D_i)).$$
This, in turn, follows from the case $p=0$ of \cite{EL}, Theorem 2,
by taking (in the notations of \cite{EL}) $L_d={\cal O}_X(D_i)$, 
$W=H^0(X, L_d)$, $N_f=K_X+D$. Indeed, $N_f=K_X+(n+1)A+n(K_X+(n+1)A)+E$,
and $n(K_X+(n+1)A)$ is nef; and $D_i$ is globally generated because
it's of the form $K_X+(n+1)A+B_i$ with $A$ very ample.

As for the surjectivity of the restriction map, one sees from the 
resolution $(\ast)$ and the Kodaira vanishing that it is enough to
check that $H^{n-1}(X,\Omega^1_X(-D))=0$, or, dually, 
$H^1(X, \Omega_X^{n-1}(D))=0$. Now the bundle $\Omega_X^{n-1}(nA)$ is
globally generated, so one can apply Griffiths' vanishing:
$H^i(X,E\otimes F)=0$ for $i>0$, where $E$ is a globally generated vector 
bundle and $F$ is a line bundle such that
 $F\otimes(det(E)\otimes K_X)^{-1}$ is ample. As
$det(\Omega_X^{n-1}(nA))=n^2A+(n-1)K_X$, the required vanishing follows.

\

Let us return to our $M\subset X$, $M=D_1\cap D_2\cap \dots \cap D_{n-2}$ with
$D_i\in |f^*{\cal O}_Q(1)|$, violating the infinitesimal Noether-Lefschetz
condition. It is clear that if $f^*{\cal O}_Q(1)$ is sufficiently ample,
then the $D_i$ shall satisfy the two conditions of our proposition, leading to
a contradiction. More precisely, suppose that $D_i$ is
of the form $2K_X+(2n+5)A$ for some very ample $A$; an easy calculation
(using the fact that $K_X+(n+1)A$ is nef) shows that the conditions 
of the proposition are satisfied provided that $n\geq 4$. 
(For $n=3$, these hold when $D=D_1=3K_X+16A$; this is already in \cite{EL}). 
So 
$D_i$ cannot
get that ample; in other words, if $m$ is such that $f^*{\cal O}_Q(1)\equiv
mL$, then $m<2k+(2n+5)r$ for $n\geq 4$ and $m<3k+16r$ for $n=3$.

\

{\bf Remark}: This can be slightly improved by remarking that $X$ is not
$\p^n$, and if $X$ is a quadric, then $f$ is an isomorphism; so one
may suppose that $X$ is neither $\p^n$ nor a quadric, in which case
already $K_X+(n-1)A$ is nef.

{}

\end{document}